\documentclass{amsart}
\usepackage{amsfonts}
\newtheorem{theorem}{Theorem}[section]

\newtheorem{corollary}[theorem]{Corollary}

\theoremstyle{definition}

\newtheorem{remark}[theorem]{Remark}

\newcommand{\ff}{\mathcal{F}}

\newcommand{\shat}{\sigma}
\newcommand{\ub}{\mathfrak{u}}
\newcommand{\ubar}{\mathfrak{u}}
\newcommand{\kbar}{\mathfrak{k}}
\newcommand{\pbar}{\mathfrak{p}}

\newcommand{\kb}{\mathfrak{k}}

\newcommand{\kp}{{\mathfrak{k}^\prime}}
\newcommand{\pp}{{\mathfrak{p}^\prime}}

\newcommand{\tbar}{\tau}

\newcommand{\Ad}{\textup{Ad}}
\newcommand{\khat}{\mathfrak{u}_+}
\newcommand{\phat}{\mathfrak{u}_-}
\newcommand{\upp}{\mathfrak{u}^{++}}
\newcommand{\upm}{\mathfrak{u}^{+-}}

\newcommand{\ump}{\mathfrak{u}^{-+}}

\newcommand{\umm}{\mathfrak{u}^{--}}

\newcommand{\real}{{\mathbb R}}

\newcommand{\hh}{\mathcal{H}}

\newcommand{\cc}{{\mathbb C}}

\newcommand{\bbar}{\left[ \begin{array}}
\newcommand{\ebar}{\end{array} \right] }
\newcommand{\bdm}{\begin{displaymath}}
\newcommand{\edm}{\end{displaymath}}
\newcommand{\beq}{\begin{equation}}
\newcommand{\beqa}{\begin{eqnarray}}
\newcommand{\beqas}{\begin{eqnarray*}}
\newcommand{\eeq}{\end{equation}}
\newcommand{\eeqa}{\end{eqnarray}}
\newcommand{\eeqas}{\end{eqnarray*}}
\newcommand{\dd}{\textup{d}}

\begin{document}

\title[Generalizations of Hilbert's  non-immersion theorem]{Results related to generalizations of Hilbert's non-immersibility theorem for the hyperbolic plane}

\author{David Brander}
\address{Department of Mathematics\\ Faculty of Science, Kobe University\\1-1 Rokkodai, Nada-ku\\ Kobe 657-8501, Japan}

\email{brander@math.kobe-u.ac.jp}

\begin{abstract}
We discuss generalizations of the well-known theorem of Hilbert that there is
no complete isometric immersion of the hyperbolic plane into Euclidean 3-space.
We show that this problem is expressed very naturally as the question of the existence
of certain homotheties of reflective submanifolds of a symmetric space.
As such,
we conclude that the only other (non-compact) cases to which this theorem could 
generalize are  the problem of isometric immersions with flat normal bundle
of the hyperbolic space $H^n$ into a Euclidean space $E^{n+k}$, $n \geq 2$,
 and the
problem of Lagrangian isometric immersions of $H^n$ into $\cc^n$, $n \geq 2$. Moreover,
there are natural compact counterparts to these problems, and for the compact
cases we prove that the theorem does in fact generalize: local embeddings exist,
but complete immersions do not.
\end{abstract}


\maketitle

\section{Introduction}

Around 1900, D Hilbert proved that there is no complete isometric
immersion of the hyperbolic plane into Euclidean 3-space \cite{hilbert1}.
Cartan studied the generalization of the problem to higher dimensions in
\cite{cartan1, cartan2}, proving that the minimal codimension needed
for even a \emph{local} isometric immersion of the hyperbolic $n$-space
$H^n$ into Euclidean space $E^{n+k}$ is $k = n-1$. 
Explicit local solutions  can be given in this codimension by a formula. 
Thus a question, which is still open for $n>2$, is whether or not \emph{complete}
isometric immersions exist $H^n \to E^{2n-1}$. For convenience, let us call this
the \emph{hyperbolic non-immersion problem}.

Building on Cartan's work, JD Moore showed
the existence of global asymptotic coordinates and the flatness of
the normal bundle for such an immersion  \cite{moore1}. The existence of these coordinates,
which had been used in the proof of Hilbert's theorem, led to the conjecture that
the higher dimensional analogue should hold.  Moreover, for space forms, these results only
depend on the fact that the extrinsic curvature (the  sectional curvature 
of the source minus that of the target) is negative. Thus, a corollary of
the existence of asymptotic coordinates, which give a covering of the immersed space 
by Euclidean space, is that, on topological grounds, there can be no
global isometric immersion of a sphere of dimension $n \geq 2$ into a sphere of 
smaller radius and of dimension $2n-1$. We call this the \emph{compact version} of the non-immersion
 problem. 
 
  In fact there really are only two versions of this negative extrinsic
 curvature non-immersion problem for space forms, because in 
 the hyperbolic version the target can equivalently be replaced by a hyperbolic
 space $H_K^{2n-1}$, with sectional curvature $K>-1$, or a sphere $S^{2n-1}(R)$ of
 any radius, $R$. The third possibility, flat immersions into a sphere, does have a
 global solution in the critical codimension, namely the Clifford torus immersion
 $E^n \to S^{2n-1}$. 
 One should point out that the \emph{positive} extrinsic curvature
 problem is not of interest, as global umbilic hypersurfaces exist.

The (hyperbolic) non-immersion  problem was studied by various people,  
such as Terng and Tenenblat \cite{terng1980, terngtenenblat},
Xavier \cite{xavier}, Pedit \cite{pedit}, 
and in several works of Y Aminov. In general, the results to date concerning
this  depend not directly on the codimension, but only
on the flatness of the normal bundle.
The compact version result mentioned above also holds in arbitrary codimension,
with this normal bundle assumption. Concerning the
hyperbolic version,  Nikolayevsky, \cite{nikolayevsky},
proved that if $M$ has constant sectional curvature $c<0$, and the fundamental
group of $M$ is non-trivial, then there is no complete isometric immersion
with flat normal bundle of $M$ into any Euclidean space. Thus, only
 the simply connected case remains.

All of the works just mentioned used the special coordinates  of Cartan and Moore
as the starting point. On the other hand, from a different
point of view, Ferus and Pedit \cite{feruspedit1996},
gave a representation of constant (non-zero) curvature submanifolds
with flat normal bundle of (non-flat) space forms as certain maps into a loop 
group, $\Lambda G$, the group of maps from the unit circle $S^1$ into a Lie group $G$.
They showed how to produce infinitely many local solutions by solving a collection
of commuting ODEs on a certain finite dimensional vector space, a standard feature
 of so-called ``integrable systems'' (as is the existence of the B\"{a}cklund 
 transformations  studied by Terng and Tenenblat).

This loop group construction was studied further in \cite{brander2} and \cite{branderrossman}.
A consequence of the loop group formulation is that it is not hard to see
that the whole construction applies much more generally, and is associated to
any pair of commuting involutions on a semisimple Lie group.  One can see
that certain questions, concerning  existence of solutions, should only depend on the
rank of a symmetric space corresponding to one of these involutions.
Thus the general context of the non-immersion problem should be determined by
identifying what is the geometric interpretation of the special submanifolds
arising in analogue to the constant curvature submanifolds of space forms.

The generalized loop group problem is investigated comprehensively in the article
\cite{reflective}. In this note, we extract the conclusions which are 
relevant to generalizations of Hilbert's theorem.  One concludes that,
within this context, the only other possible
generalization, besides that of isometric immersions with flat normal bundle
of $H^n$ into $E^{n+k}$, is the problem of isometric Lagrangian immersions
of $H^n$ into $\cc^n$. For these cases, local solutions can be constructed
by integrable systems methods, and the author
does not know whether global solutions exist. 
In all other potential cases (described below) local solutions do not exist. 

We also determine all possible generalizations of the compact version of the
problem which are: the problem of isometric immersions with flat normal bundle of a sphere 
$S^n(R)$, of radius $R>1$, into a unit sphere $S^{n+k}$, $k \geq n-1$,
 and
the problem of isometric Lagrangian immersions of $S^n(R)$, $R>1$,
  into $\cc P^n$.  We prove that these have local, but no global, solutions.
  
\begin{remark} A stronger version of Hilbert's theorem was proved by 
Efimov in \cite{efimov1964}. He proved that there is no complete 
isometric immersion in
$E^3$ of a surface whose Gauss curvature is bounded above by some negative
constant.  Generalizations to higher dimensions of this stronger result 
have been in the direction of hypersurfaces \cite{xaviersmyth}, rather than
to codimension $n-1$.
\end{remark}

\section{Generalizations of the Compact Version}
\subsection{The Loop Group Construction}
We primarily describe the compact case in this article. 
Full details of all cases can be found in \cite{reflective}.
The basis of the method is that given an immersion into a homogeneous 
space, $f: M \to G/H$, one can lift $f$ to a frame $F: M \to G$. Up to an
irrelevant isometry of $G/H$, the immersion
$f$ is completely determined by the pull-back to $M$ of the Maurer-Cartan form
of $G$, denoted by $F^{-1} \dd F$, and called the Maurer-Cartan form of $F$.
One can study special submanifolds by choosing an appropriately adapted frame,
and the geometry is encoded in the Maurer-Cartan form.
 
Let $G$ be a complex semisimple Lie group, and $U$ a real form
defined as the fixed point subgroup of a complex antilinear involution $\rho$ of $G$.
Let $\sigma$ and $\tau$ be a pair of involutions of $G$, and suppose that all three
involutions commute.  Let $\Lambda G$ denote the group of maps from the unit circle
$S^1$ to $G$, of a suitable class, so that $\Lambda G$ is a Banach Lie group.
We denote the $S^1$ parameter by $\lambda$.  Extend the involutions to $\Lambda G$ by
the formulae:
\beqas
(\rho x)(\lambda) := \rho(x(\bar{\lambda})), \hspace{1cm}
(\tau x)(\lambda) := \tau(x(-\lambda^{-1})), \hspace{1cm}
(\sigma x)(\lambda) := \sigma(x(-\lambda)),
\eeqas
where $x:S^1 \to G$ is any element of $\Lambda G$.
Now define a subgroup $\hh$ of $\Lambda G$ as the set of elements which are fixed by all
three involutions,
\bdm
\hh := \{ x \in \Lambda G ~|~ \rho x = \tau x = \sigma x = x \}.
\edm
The point of this loop group construction will be that the maps we consider
give families of specially adapted frames, $F_\lambda$, whose Maurer-Cartan forms 
have a particular expression (see (\ref{alphamcf}) below).  The description in
terms of the extended involutions above is  important for the problem 
of \emph{constructing} solutions, where this description fits naturally into very 
general methods, which we shall not describe.

The Lie algebra, $Lie(\hh)$, of $\hh$ consists of Laurent polynomials in $\lambda$, 
$\sum X_i \lambda ^i$, with
coefficients $X_i$ in certain subspaces of the Lie algebra $\mathfrak{g}$ of $G$, and an appropriate convergence
condition.  Let $\hh^0 := \hh \cap G$, the subgroup of constant loops. The type
of loop group maps we consider are smooth maps $f: M \to \hh/\hh^0$, where $M$ is a 
simply connected manifold, and $\hh/\hh^0$ is
the left coset homogeneous space. Moreover, we impose the restriction that
for any lift, $F: M \to \hh$, of $f$,  the Maurer-Cartan form  $F^{-1} \dd F$ 
of $F$ is a Laurent polynomial (the coefficients of which are $\mathfrak{g}$-valued 1-forms)
whose highest and lowest powers of $\lambda$ are
1 and $-1$ respectively.  Denote the set of such maps by:
\bdm
\ff(M) := \{f: M \to \hh/\hh^0 ~|~ F^{-1} \dd F = \sum_{i=-1}^1 \alpha_i \lambda^i, ~\forall ~\textup{lifts}~F \}.
\edm
Clearly, if we fix a value of the loop parameter, $\lambda$, an element
$f \in \ff(M)$, gives a map $f_\lambda: M \to G/\hh^0$, with 
corresponding frames $F_\lambda$.
 The restriction on the Maurer-Cartan
form of $F$ ensures that its dependence on $\lambda$ extends 
holomorphically to the punctured complex plane $\cc^*$, so we can  consider $f_\lambda$
 for non-zero real values $\lambda \in \real^*$. 
The condition $\rho F = F$ implies that, for real values of $\lambda$, 
$F_\lambda$ takes its values in the real form $U$. It is  clear 
that $\hh^0$ is the fixed
point subgroup of $G$ with respect to the three involutions $\rho$, $\tau$ and $\sigma$.
In fact $\hh^0 = K \cap U_+$, where
\bdm
K = U_\tau, \hspace{1cm} U_+ = U_\sigma,
\edm
are the fixed point subgroups  with respect to $\tau$ and $\sigma$.

Thus, for $\lambda \in \real^*$, $f_\lambda$ is a map 
\bdm
M \to \frac{U}{K \cap U_+},
\edm
and one can consider projections to either of the symmetric spaces $U/K$ or  $U/U_+$.
Denote these projections by $\bar{f}_\lambda : M \to U/K$ and 
$\hat{f}_\lambda: M \to U/U_+$.

Set 
\bdm
R_\lambda = \Big | \frac{\lambda + \lambda^{-1}}{2} \Big |.
\edm
Note that $R_\lambda \geq 1$ for $\lambda \in \real^*$. 

\begin{theorem}
\cite{feruspedit1996, reflective} \label{thm1} 
 Projection to $U/K$:
Suppose that $M$ has dimension
$n$, and let $f \in \ff(M)$. Suppose that
 the projection of the map obtained at $\lambda =1$,
$\bar{f}_1: M \to U/K$ is regular (i.e. an immersion).
 Then so is $\bar{f}_\lambda$ for any 
other value of $\lambda \in \real^*$, and
\begin{enumerate}
\item [(i)] \label{item1}
Suppose $U = SO(n+k+1)$, $\sigma = \Ad_P$,  $\tau = \Ad_Q$,  where
\bdm
 P = \bbar{cc} I_n & 0 \\ 0 & -I_{k+1} \ebar,
 \hspace{1cm} Q = \bbar{cc} I_{n+k} & 0 \\ 0 & -1 \ebar
\edm
and $I_l$ denotes an $l \times l$ identity matrix.

Then $U/K = S^{n+k}$, and
$\bar{f}_\lambda: M \to S^{n+k}$, with the induced metric, 
 is an isometric immersion with flat normal bundle of a part of a sphere 
$S^n(R_\lambda)$, of radius $R_\lambda$.  
\item [(ii)] \label{item2}
Suppose $U = SU(n+1)$, represented by the matrix subgroup of $SO(2n+2)$ consisting
of all matrices of the form $\tiny{\bbar {cc} A & -B\\ B & A \ebar}$, where $A$ and
$B$ are $(n+1) \times (n+1)$, and such that $\det(A+iB) = 1$. 
Let $\shat = \Ad _P$, for
 $P ={\textup{diag}(I_{n+1},-I_{n+1})}$, and $\tbar = \Ad_Q$, for 
 $Q = \textup{diag}(I_n, -1, I_n, -1)$. 

Then $U/K = \cc P^n$, and $\bar{f}_\lambda : M \to \cc P^n$, with the induced metric, 
 is a Lagrangian  isometric immersion 
 of a part of a sphere 
$S^n(R_\lambda)$, of radius $R_\lambda$.
\end{enumerate}
Conversely, in both cases, any such isometric immersion, with $R>1$,
can be represented by such an element $f \in \ff(M)$.
\end{theorem}
Note that for the converse one needs $R$ to be strictly greater than 1,
for reasons which will become clear in the proof outlined below.

\begin{theorem} \cite{reflective} \label{thm2} 
{Projection to $U/U_+$:}
Let $f \in \ff(M)$.  For $\lambda \in \real^*$, let 
$\hat{f}_\lambda :M \to U/U_+$ be the map obtained by projection. In the limit
as $\lambda \to \infty$, $\hat{f}_\lambda$ is asymptotic to a flat:
that is, a flat totally geodesic submanifold of the  symmetric
space $U/U_+$.  If the projection $\bar{f}_\lambda : M \to U/K$
is an immersion, then so is $\hat{f}_\lambda$. In this case,
 $M$ admits a flat metric.
\end{theorem}

Note that in \cite{reflective} (Proposition 4.2) it is stated, too strongly,
 that the projection
$\hat{f}_\lambda$ is a curved flat for all real $\lambda$. In fact this is only
true  in the limit as $\lambda$ approaches $0$ or $\infty$.

A flat is obtained by exponentiating an Abelian subalgebra of $\mathfrak{u}_-$.
If the symmetric space $U/U_+$ is Riemannian, that is, if $U_+$ is compact, then
 the dimension of such a  subalgebra  is, by definition,
  no greater than
 the rank of $U/U_+$. Hence, Theorem \ref{thm2} has the following 
\begin{corollary} \label{cor1}
Let $f \in \ff(M)$, and suppose that the map $\bar{f}_\lambda: M \to U/K$ is regular.
Then, if $U_+$ is compact, the dimension of $M$ is no greater than the rank 
of $U/U_+$.
\end{corollary}

Applying this condition to the first case in Theorem \ref{thm1}, we conclude
that the minimal codimension needed for an isometric immersion $S^n(R) \to S^{n+k}$,
where $R>1$, is $k = n-1$. Thus Corollary \ref{cor1}
 gives a natural explanation for this  known result.

For the second case in
 Theorem \ref{thm1}, the symmetric space $U/U_+$ is just
$SU(n+1)/SO(n+1)$, which has rank $n$, thus local solutions are not ruled out.
In fact Corollary \ref{cor1} has a converse:
\begin{theorem} \label{thm3} \cite{reflective}
If $\textup{Dim}(M) \leq \textup{Rank}(U/U_+)$, 
then solutions $f \in \ff(M)$ can
always be constructed, at least locally, which have regular projections to $U/K$.
\end{theorem}
Hence we conclude that local 
Lagrangian isometric immersions $S(R)^n \to \cc P^n$ do indeed exist for $R>1$.
 
Finally, the flat metric of Theorem \ref{thm2} implies the existence  
of  a topological covering by a Euclidean space. We therefore conclude:
\begin{corollary} \label{cor2}
If the dimension of $M$ is greater than 1, then none of the solutions $\bar{f}_\lambda$
in either case of Theorem \ref{thm1}  can be complete.
\end{corollary}

\subsection{The Proofs of Theorems \ref{thm1} and \ref{thm2}}
Let
\bdm
\ubar = \kbar \oplus \pbar = \khat \oplus \phat
\edm
be the canonical decompositions of the Lie algebra of $U$, 
associated to $\tau$ and $\sigma$ respectively.
We have the orthogonal (with respect to the Killing-form of $\mathfrak{u}$) decomposition
\bdm
\ub = \upp \oplus \upm \oplus \ump \oplus \umm, 
\edm
where
\beqas
\ub^{++} =   \kb \cap \khat =: \kp, && \hspace{1cm}
 \ub^{+-} =  \kb \cap \phat,\\
 \ub^{--} =  \pbar \cap \phat =: \pp, && \hspace{1cm}
 \ub^{-+} = \pbar \cap \khat =: \mathfrak{p}^{\prime \perp}.
 \eeqas
The starting point for the proof of both theorems is the  fact
that if $F: M \to \hh$ is a lift of an element $f \in \ff(M)$, then, for $\lambda \in \real^*$, the 
Maurer-Cartan form, $\alpha^\lambda := F^{-1}\dd F$,
 of $F$ has the following expansion:
\beq  \label{alphamcf}
\alpha^\lambda= \alpha_0^{++} + \alpha_1^{+-} (\lambda - \lambda^{-1}) +  \alpha_1^{--}(\lambda + \lambda^{-1}),
\eeq
where
\bdm
\nonumber \alpha_0^{++} \in \ub^{++} \otimes \Omega(M), \hspace{.5cm}
 \alpha_1^{+-} \in   \ub^{+-} \otimes \Omega(M),   \hspace{.5cm}
  \alpha_1^{--} \in \ub^{--}\otimes \Omega(M),
\edm  
and $\Omega(M)$ denotes the space of 1-forms on $M$.
The expansion (\ref{alphamcf}) is deduced from the fact that $\alpha^\lambda$ is
fixed by the infinitesimal versions of $\sigma$ and $\tau$.

\subsubsection{Proof of Theorem \ref{thm1}}
Recall that $F_\lambda$ is $U$-valued for $\lambda \in \real^*$.
We can use the frame $F_\lambda$ to
give vector bundle isomorphisms between the tangent and normal bundles for $\bar{f}_\lambda$,
$\lambda \neq 1$, and those at $\lambda =1$, by left translating them 
to $\mathfrak{p}$ via $F_1^{-1}$ and $F_\lambda^{-1}$ respectively as follows:
the tangent space to $U/K$ at $f_\lambda(x)$ is identified with $\mathfrak{p}$
by left multiplication by $F_\lambda(x)^{-1}$, in the standard way.
Giving the immersions the metric induced from the 
standard Killing form metric on $U/K$, 
the coframe is thus given by the $\mathfrak{p}$ component of this 1-form, namely
$\alpha_1^{--}(\lambda + \lambda^{-1})$, from which one deduces that the induced
metric for $\bar{f}_\lambda$ 
is just a scalar multiple of the metric at $\lambda =1$.
In both cases of Theorem \ref{thm1},
$\textup{Dim}(M) = \textup{Dim}(\mathfrak{p}^\prime$), and so it follows 
that the tangent space to $f_\lambda(M)$ and the normal space
are given respectively by $\mathfrak{p}^\prime$ and $\mathfrak{p}^{\prime \perp}$
under the above identification.
The connection 1-form for $\bar{f}_\lambda$ is given by the projection to $\mathfrak{k}$
of $\alpha^\lambda$, and this splits into the connections on the tangent and normal 
bundles as well as the second fundamental form. 
One has the relations
\beqas  \label{1strelation}
~[\kb \cap \khat, \pp] \subset \pp, \hspace{1cm}
[\kb \cap \khat, \pp^\perp] \subset \pp^\perp,\\
~[\kb \cap \phat, \pp] \subset \pp^\perp, \hspace{1cm} 
[\kb \cap \phat, \pp^\perp] \subset \pp, \label{2ndrelation}
\eeqas
from which it is not
difficult to see that the second fundamental form is given by the 
$\mathfrak{k} \cap \mathfrak{u_-}$ component, namely
$\alpha_1^{+-} (\lambda - \lambda^{-1})$, which implies that, at $\lambda =1$,
we have (part of) a totally geodesic submanifold $N$ of $U/K$. Since
 the coframe takes values in
$\mathfrak{p}^\prime$, $N$ must be the projection
to $U/K$ of $\exp(\mathfrak{p}^\prime)$. In the first case,
 $N$ is a totally geodesic sphere $S^n$ in $SO(n+k+1)/SO(n+k) = S^{n+k}$.
 In the second case, $N$ is a totally geodesic 
Lagrangian submanifold of $\cc P^n$.

Finally, the 1-form $\alpha_0^{++}$ is the sum of the tangential and normal
connections of $\bar{f}_\lambda$, which does not depend on $\lambda$. 
This means that, under the isomorphisms between the respective tangent and
normal bundles for different values of $\lambda$ given above, the connections
on these bundles are also preserved. 
  For the first case, this means the normal bundle is flat, as the totally geodesic
  sphere $N$ has this property. For the second case, where $N$ is Lagrangian,
it follows that, for $\lambda \neq 1$, $\bar{f}_\lambda$ is also Lagrangian, as the complex
structure is given on $\mathfrak{p}$, and the tangent and normal bundles for
different values of $\lambda$ are identified in $\mathfrak{p}$.

The converse can be obtained in a fairly
straightforward manner, by choosing appropriately adapted frames for the type
of immersions required. Note that one cannot have $R=1$ for the converse,
because this corresponds to $\lambda - \lambda^{-1}=0$ and we would not
know what 1-form to insert for $\alpha_1^{+-}$ in (\ref{alphamcf}).

\subsubsection{Proof of Theorem \ref{thm2}}
For the projection, $\hat{f}$, to $U/U_+$, the coframe is
given by the  $\mathfrak{u}_-$ component of $\alpha_\lambda$, namely
$\beta =  \alpha_1^{+-} (\lambda - \lambda^{-1}) +  \alpha_1^{--}(\lambda + \lambda^{-1})$.
The limiting coframe as $\lambda \to \infty$, is proportional to the 1-form
\bdm
\beta_\infty = \alpha_1^{+-} + \alpha_1^{--}.
\edm 

Now for any value of $\lambda$, $\alpha^\lambda$ must satisfy the Maurer-Cartan equation, 
$\dd \alpha + \alpha \wedge \alpha = 0$, this condition being  equivalent to the existence
of a map $F_\lambda$ such that $\alpha^\lambda = F_\lambda^{-1} \dd F_\lambda$. 
The fact that this holds for \emph{all} values of $\lambda$ implies
the curved flat
equation, $\beta_\infty \wedge \beta_\infty =0$, 
or equivalently, that the matrix components of
$\beta$ all commute, from which one can 
deduce that $\hat{f}_\lambda$ is asymptotic to a flat as $\lambda \to \infty$.
The condition that $\bar{f}_\lambda$ be regular is just that its coframe, $\alpha_1^{--}$, consists of
$n$ linearly independent 1-forms.  The condition that $\hat{f}_\lambda$ be regular
is that $\beta$ has the same property, and this clearly follows from the regularity of
$\bar{f}_\lambda$. 

 Finally, one can show that the metric given by
$(X,Y) := \langle \beta_\infty((F_\lambda)_* X), \beta_\infty((F_\lambda)_* Y) \rangle$,
where $\langle, \rangle$ is the Killing metric on $\mathfrak{u}_-$, is a flat
metric on $M$. One way to see this is that, as shown in \cite{reflective}, there
is, locally associated to $F_\lambda$, a curved flat $f_+ : M \to U/U_+$ 
(see the proof of Proposition 5.2), and the coframe of $f_+$ is given by
$\psi = \textup{Ad}_C \beta_\infty$, where $C$ takes values in $U_+$. Since
$U_+$ acts by isometries on $\mathfrak{u}_-$, the metric $(,)$ is the same
as that induced by the curved flat. It is shown in \cite{feruspedit1996II}
that such a metric is flat.

\subsection{Generalizations to Other Symmetric Spaces}
 It is known that such a pair of involutions,
$\tau$ and $\sigma$ define a \emph{reflective submanifold} $N$ of the symmetric
space $U/K$; that is, a totally geodesic submanifold which has an external symmetry.
This is just the projection to $U/K$ of $\exp(\mathfrak{p}^\prime)$, mentioned
in the previous section.  In fact all connected totally geodesic submanifolds
of symmetric spaces are given by $\exp(\mathfrak{p}^\prime)$ for some 
vector subspace $\mathfrak{p}^\prime$ of $\mathfrak{p}$ which is closed under
the Lie triple product, 
$[\mathfrak{p}^\prime,[ \mathfrak{p}^\prime, \mathfrak{p}^\prime]] \subset
 \mathfrak{p}^\prime$. A reflective submanifold has the additional 
 property that the orthogonal complement in $\mathfrak{p}$,
  denoted $\mathfrak{p}^{\prime \perp}$,
  is also closed under the triple product. This is equivalent
 to the existence of the second involution $\sigma$.
 
  The argument outlined above applies in all cases,
assuming that $\textup{Dim}(M) = \textup{Dim}(\mathfrak{p}^\prime)$. In general,
the projections
$\bar{f}_\lambda$ of elements $f \in \ff(M)$ correspond to 
certain homotheties of the
reflective submanifold $\bar{f}_1$ obtained at $\lambda =1$, keeping the normal
bundle isomorphic.  We have shown that, in the compact case, globally, there is
no such homothety for any reflective submanifold.  To check whether local solutions
exist for other cases is just a matter of comparing the dimension of the  reflective
submanifold with the rank of the associated
second symmetric space, $U/U_+$.
Now reflective submanifolds of symmetric spaces 
 were studied and classified by DSP Leung in \cite{leung1974, leung1975,leung1979},
 and there are many cases.  However, it turns out that in all the other cases
 the rank is too small for local solutions to exist. Hence we conclude that
Corollary \ref{cor2} contains all possible generalizations (to reflective submanifolds
of simply connected, compact, irreducible, Riemannian symmetric spaces)
of the compact version of the Hilbert theorem.

\section{The Hyperbolic Case}
For the hyperbolic case, the problem we are interested in, namely negative extrinsic
curvature, corresponds to homotheties of the reflective submanifold by a factor
$R<1$, rather than greater than 1.  This problem also has a loop group formulation,
which differs from that described above only in that the loops are real-valued for
values of the parameter $\lambda$ in $S^1$, rather than $\real^*$.  However, by
evaluating such a loop group map for values of $\lambda$ in $\real^*$,
instead of in $S^1$,
one obtains a similar situation to that of the compact case (although for a
different, non-Riemannian, symmetric space $\tilde{U}/\tilde{K}$) and one
can obtain analogous results to those described in the compact case, with the
exception of the non-existence of global solutions, which remains an open problem.
In particular, it is shown in \cite{reflective} that the only cases (of reflective
submanifolds of simply connected, non-compact, irreducible, Riemannian symmetric
spaces) where local solutions exist are:
\begin{enumerate}
\item [(i)] 
$U/K = H^{n+k}$, and
$\bar{f}_\lambda: M \to H^{n+k}$, $k \geq n-1$, with the induced metric, 
 is an isometric immersion with flat normal bundle of a part of a hyperbolic space 
$H^n_c$ of constant sectional curvature $c<-1$.  
\item [(ii)] 
$U/K = \cc H^n$, and $\bar{f}_\lambda : M \to \cc H^n$, with the induced metric, 
 is a Lagrangian  isometric immersion 
 of a part of a hyperbolic space 
$H^n_c$ of constant sectional curvature $c<-1$.  
\end{enumerate}
In both cases local solutions exist and can be constructed using integrable
systems methods.

Equivalently, one can replace the target spaces in cases (i) and
(ii) with, respectively, the Euclidean space, $E^{n+k}$, and complex Euclidean space $\cc^n$,
by an argument which was given in \cite{branderrossman}. That article dealt
with the first case only, but the proof is easily adapted to the second case.
It essentially involves dilating the target space, while keeping the metric on
the immersed space constant, until, in the limit, an immersion into flat space
is obtained.

\providecommand{\bysame}{\leavevmode\hbox to3em{\hrulefill}\thinspace}
\providecommand{\MR}{\relax\ifhmode\unskip\space\fi MR }
\providecommand{\MRhref}[2]{%
  \href{http://www.ams.org/mathscinet-getitem?mr=#1}{#2}
}
\providecommand{\href}[2]{#2}

\end{document}